\documentclass[epj]{webofc}
\usepackage[varg]{txfonts}   
\woctitle{Mathematical Modeling and Computational Physics 2015}

\def\md{\mathrm{d}}

\begin{document}
\title{Numerical Solution of a Nonlinear Integro-Differential Equation}
%
%

\author{J\'an Bu\v{s}a\inst{1}\fnsep\thanks{\email{jan.busa@tuke.sk}} 
\and
        Michal Hnati\v{c}\inst{2,3,4}\fnsep\thanks{\email{hnatic@saske.sk}} 
\and
        Juha Honkonen\inst{5}\fnsep\thanks{\email{juha.honkonen@helsinki.fi}}
\and
 Tom\'a\v{s} Lu\v{c}ivjansk\'y\inst{2,6}
}

\institute{Department of Mathematics and Theoretical Informatics, FEE\&I, Technical University, 
	Ko\v{s}ice, Slovakia 
\and
         Faculty of Sciences, P.~J.~\v{S}afarik University, Ko\v{s}ice, Slovakia
         \and
         Institute of Experimental Physics SAS, Ko\v{s}ice, Slovakia
         \and
         Bogoliubov Laboratory of Theoretical Physics, JINR,
         141980 Dubna, Moscow Region, Russia 
\and
    Department of Military Technology, National Defence University,  Helsinki, Finland
    \and
Fakult\"at f\"ur Physik, Universit\"at Duisburg-Essen, D-47048  
    Duisburg, Germany
       }

\abstract{%
An algorithm for the numerical solution of a nonlinear integro-differential equation arising in the single-species annihilation reaction $A+A \to \emptyset$ modeling is discussed. Finite difference method together with the linear approximation of the unknown function is considered. For divergent integrals presented in the equation for dimension $d=2$ a regularization is used. Some numerical results are presented.}
\maketitle
\section{Introduction}
\label{intro}

The irreversible annihilation reaction $A + A \to \emptyset$ is a fundamental model of non-equilibrium physics. The reacting $A$ particles are assumed to perform chaotic motion due to diffusion or some external advection field such as atmospheric eddy \cite{HHL}. Many reactions of this type are observed in diverse
chemical, biological or physical systems \cite{r1,r2}.

In \cite{HHL} the advection of reactive scalar using random velocity field
generated by the stochastic Navier-Stokes equation, which
is used for production a velocity field corresponding to
thermal fluctuations \cite{r15,r16} and a turbulent velocity field
with the Kolmogorov scaling behavior \cite{r17} is studied by three
of the present authors, and the integro-differential equation for the number density is derived. No influence of the reactant on the velocity field itself is assumed.

In this paper we present some initial experiences with the numerical solution of the integro-differential equation mentioned above.

\section{Problem Formulation}
\label{sec-1}

In \cite{HHL} an integro-differential equation (72) for the mean number density $a(t)$ of chemically active molecules in anomalous kinetics of single-species annihilation reaction $A+A \to \emptyset$  
\begin{equation}\label{eqHHL72}
\frac{\md a(t)}{\md t}=-2\lambda u\nu \mu^{-2\Delta} Z_4\, a^2(t)+4\lambda^2 u^2 \nu^2 \mu^{-4\Delta} \int_0^t \frac{a^2(t')\,\md t'}{\big[8\pi u\nu(t-t')\big]^{d/2}}, \qquad a(0)=a_0
\end{equation}
is derived. The integral in Eq.~\eqref{eqHHL72} diverges at the upper limit $t'$ in space dimensions $d \geq 2$.

We will study a numerical solution to the Initial Value Problem for the equation:
\begin{equation}\label{eq:ide}
\frac{\md a(t)}{\md t}=-2\lambda D a^2(t)+4\lambda^2 D^2 \int_0^t \frac{a^2(t')\,\md t'}{\big[8\pi D(t-t')\big]^{d/2}}, \qquad a(0)=a_0,
\end{equation}
which corresponds to $D=u \nu$, $\Delta=0$, and $Z_4=1$.

\section{Case $d=2$}
For $d=2$ the singularity in the integral on the right side of Eq.~\eqref{eq:ide} at $t'=t$ is divergent. This is a consequence of the UV divergences in the model above the critical dimension $d_c = 2$, and near the critical dimension is remedied by the UV renormalization of the model \cite{HHL}. In this paper we use another approach to overcome this problem. We will use the following \emph{regularization}:
\begin{equation}\label{eq:ider2}
\frac{\md a}{\md t}=-2\lambda D a^2+4\lambda^2 D^2 \int_0^t \frac{a^2(t')\,\md t'}{8\pi \{D(t-t')+\ell^2\}}, \qquad a(0)=a_0.
\end{equation}
We rewrite Eq.~\eqref{eq:ider2} to the form 
\begin{equation}\label{eq:ideabg2}
\frac{\md a}{\md t}=-\alpha a^2+\beta \int_0^t \frac{a^2(t')\,\md t'}{t-t'+\gamma}
\end{equation}
with
\[
\alpha=2\lambda D, \qquad \beta=\frac{\alpha^2}{8\pi D}, \qquad \gamma=\frac{\ell^2}{D}.
\]

Let us consider the Eq.~\eqref{eq:ideabg2}. We will use the difference method to solve it numerically. We will consider time discretization with the time step $\Delta t$:
\begin{equation}\label{eq:dt}
t_k= k\cdot \Delta t, \quad a(t_k)=a_k, \qquad k=0, 1, 2, \dots
\end{equation}
For the discretization of the left side we will use two formula -- the first and the second order finite difference  approximations:
\begin{equation}\label{eq:dadt}
\left.\frac{\md a}{\md t}\right|_{t=t_k}\approx \dfrac{a_k-a_{k-1}}{\Delta t}\qquad\mbox{or}\qquad \left.\frac{\md a}{\md t}\right|_{t=t_k}\approx \dfrac{3a_k-4a_{k-1}+a_{k-2}}{2\Delta t}.
\end{equation}
For the right-side integral approximation we will use piecewise linear approximation of the function $a(u)$:
\begin{equation}\label{eq:int}
\int_0^{t_k} \frac{a^2(t')\,\md t'}{t_k-u+\gamma}=\sum_{i=0}^{k-1}\int_{i\cdot\Delta t}^{(i+1)\cdot\Delta t} \frac{a^2(t')\,\md t'}{t_k-t'+\gamma}\approx
\end{equation} 
\[
\approx\sum_{i=0}^{k-1}\int_{i\cdot\Delta t}^{(i+1)\cdot\Delta t}\frac{\left[t'\cdot(a_{i+1}-a_i)/\Delta t+(i+1)\cdot a_i-i\cdot a_{i+1}\right]^2\,\md t'}{t_k-t'+\gamma}.
\]

Integrals in the right side of Eq.~\eqref{eq:int} we calculate analytically using
\begin{equation}\label{eq:gd}
\gamma=\delta\cdot \Delta t
\end{equation}
in the following way
\[
\int_{i\cdot\Delta t}^{(i+1)\cdot\Delta t} \frac{\left[t'\cdot(a_{i+1}-a_i)/\Delta t+(i+1)\cdot a_i-i\cdot a_{i+1}\right]^2\,\md t'}{t_k-t'+\gamma}=
\]
{
	\[
	=\int_{i\cdot\Delta t}^{(i+1)\cdot\Delta t} \frac{\left[(t'-k\cdot \Delta t-\delta\cdot \Delta t)\cdot(a_{i+1}-a_i)/\Delta t+(k+\delta)\cdot(a_{i+1}-a_i)+(i+1)\cdot a_i-i\cdot a_{i+1}\right]^2\,\md t'}{k\cdot\Delta t-t'+\delta\cdot\Delta t}=
	\]}
{
	\[
	=\int_{i\cdot\Delta t}^{(i+1)\cdot\Delta t}\frac{\left[(k+\delta-i)\cdot a_{i+1}-(k+\delta-i-1)a_i-(k\cdot \Delta t-t'+\delta\cdot\Delta t)\cdot(a_{i+1}-a_i)/\Delta t\right]^2\,\md t'}{k\cdot\Delta t-t'+\delta\cdot\Delta t}=
	\]}
\[
=\left[(k+\delta-i) a_{i+1}-(k+\delta-i-1)a_i\right]^2\cdot\int_{i\cdot\Delta t}^{(i+1)\cdot\Delta t}\!\!\! \frac{\md t'}{k\cdot\Delta t-t'+\delta\cdot\Delta t}-
\]
\[
-\frac{2\left[(k+\delta-i) a_{i+1}-(k+\delta-i-1)a_i\right]\cdot(a_{i+1}-a_i)}{\Delta t}\int_{i\cdot\Delta t}^{(i+1)\cdot\Delta t}\frac{(k\cdot \Delta t-t'+\delta\cdot\Delta t)\,\md t'}{(k\cdot\Delta t-t'+\delta\cdot\Delta t)}+
\]
\[
+\frac{(a_{i+1}-a_i)^2}{(\Delta t)^2}\int_{i\cdot\Delta t}^{(i+1)\cdot\Delta t} \frac{(k\cdot \Delta t-t'+\delta\cdot\Delta t)^2\,\md t'}{(k\cdot\Delta t-t'+\delta\cdot\Delta t)}=
\]
\[
=\left[(k+\delta-i)a_{i+1}-(k+\delta-i-1)a_i\right]^2\cdot\left[-\ln(k\cdot\Delta t- t'+\delta\cdot\Delta t)\right]_{i\cdot\Delta t}^{(i+1)\cdot\Delta t}-
\]
\[
-\frac{2\left[(k+\delta-i) a_{i+1}-(k+\delta-i-1)a_i\right](a_{i+1}-a_i)}{\Delta t}\cdot \Delta t
+\frac{(a_{i+1}-a_i)^2}{(\Delta t)^2}\cdot\left[\frac{(k\cdot \Delta t-t'+\delta\cdot\Delta t)^2}{-2}\right]_{i\cdot\Delta t}^{(i+1)\cdot\Delta t}=
\]
\begin{equation}\label{eq:intdr2}
=\Big\{\big[(k+\delta-i)a_{i+1}-(k+\delta-i-1)a_i\big]^2\cdot\ln\frac{k+\delta-i}{k+\delta-i-1}-\Big(k+\delta-i+\frac{1}{2}\Big)(a_{i+1}-a_i)^2-2a_i(a_{i+1}-a_i)\Big\},
\end{equation}
$i=0$, 1, \dots, $k-1$.

Putting approximations Eq.~\eqref{eq:dadt}, Eq.~\eqref{eq:int}, and Eq.~\eqref{eq:intdr2} into Eq.~\eqref{eq:ideabg2} we arrive at the quadratic equations Eq.~\eqref{qeqk} and Eq.~\eqref{qeq0} in more explicit form with respect to $a_k$:
\begin{equation}\label{eq:quadl}
0=\dfrac{3a_k-4a_{k-1}+a_{k-2}}{2\Delta t}+\alpha\cdot a_k^2-
\end{equation}
\[
-\beta\cdot\sum_{i=0}^{k-2}\Big\{\big[(k+\delta-i)a_{i+1}-(k+\delta-i-1)a_i\big]^2\cdot\ln\frac{k+\delta-i}{k+\delta-i-1}-\Big(k+\delta-i+\frac{1}{2}\Big)(a_{i+1}-a_i)^2-2a_i(a_{i+1}-a_i)\Big\}-
\]
\[
-\beta\cdot\Big\{\big[(1+\delta)a_{k}-\delta a_{k-1}\big]^2\cdot\ln\frac{1+\delta}{\delta}-\Big(1+\delta+\frac{1}{2}\Big)(a_{k}-a_{k-1})^2-2a_{k-1}(a_{k}-a_{k-1})\Big\}
\]
or in the standard form
\[
\bigg[\alpha-\beta\cdot\Big\{(1+\delta)^2\cdot\ln\frac{1+\delta}{\delta}
-\Big(\frac{3}{2}+\delta\Big)\Big\}\bigg]\cdot a_k^2+
\bigg[\dfrac{3}{2\Delta t}-\beta\cdot a_{k-1}\cdot\Big\{2\delta+1-2(1+\delta)\delta \cdot\ln\frac{1+\delta}{\delta}
\Big)\Big\}\bigg]\cdot a_k-
\]
\begin{equation}\label{qeqk}
-\dfrac{4a_{k-1}-a_{k-2}}{2\Delta t}-\beta\cdot a_{k-1}^2\cdot\Big\{\frac{1}{2}-\delta+\delta^2\cdot\ln\frac{1+\delta}{\delta}
\Big\}-
\end{equation}
\[
-\beta\cdot\sum_{i=0}^{k-2}\Big\{\big[(k+\delta-i)(a_{i+1}-a_i)+a_i\big]^2\cdot\ln\frac{k+\delta-i}{k+\delta-i-1}
-(a_{i+1}-a_i)\Big[(k+\delta-i)(a_{i+1}-a_i)+a_i+\frac{a_{i+1}+a_i}{2}\Big]\Big\}
=0.
\]

For the first step we have
\[
\frac{a_1-a_0}{\Delta t}+\alpha a_1^2-\beta \Big\{\big[(1+\delta)a_{1}-\delta a_0\big]^2\cdot\ln\frac{1+\delta}{\delta}-\Big(1+\delta+\frac{1}{2}\Big)(a_{1}-a_0)^2-2a_0(a_{1}-a_0)\Big\}=0
\]
\newpage
\noindent 
or in the standard form
\[
\bigg[\alpha-\beta\cdot\Big\{(1+\delta)^2\cdot\ln\frac{1+\delta}{\delta}
-\Big(\frac{3}{2}+\delta\Big)\Big\}\bigg]\cdot a_1^2
+\bigg[\dfrac{1}{\Delta t}-\beta\cdot a_{0}\cdot\Big\{2\delta+1-2(1+\delta)\delta \cdot\ln\frac{1+\delta}{\delta}
\Big)\Big\}\bigg]\cdot a_1-
\]
\begin{equation}\label{qeq0}
-\dfrac{a_{0}}{\Delta t}-\beta\cdot a_{0}^2\cdot\Big\{\frac{1}{2}-\delta+\delta^2\cdot\ln\frac{1+\delta}{\delta}
\Big\}=0.
\end{equation}

The algorithm consists of the successive calculation of the values $a_k$, $k=1,$ 2, \dots, $K$ solving the quadratic equation Eq.~\eqref{qeq0} at the first step using the value $a_0=a(0)$, and further  solving Eqs.~\eqref{qeqk} using previously determined values $a_i$, $i=0$, \dots, $k-1$.  

\section{Numerical results}

Below the results for $a(0)=5000$, $\lambda=0.1$, $D=0.4$ and $\ell=0.001$ are presented. All calculations are done with the uniform time step $\Delta t$.

\begin{figure}
	\centering
	\includegraphics[width=1\textwidth]{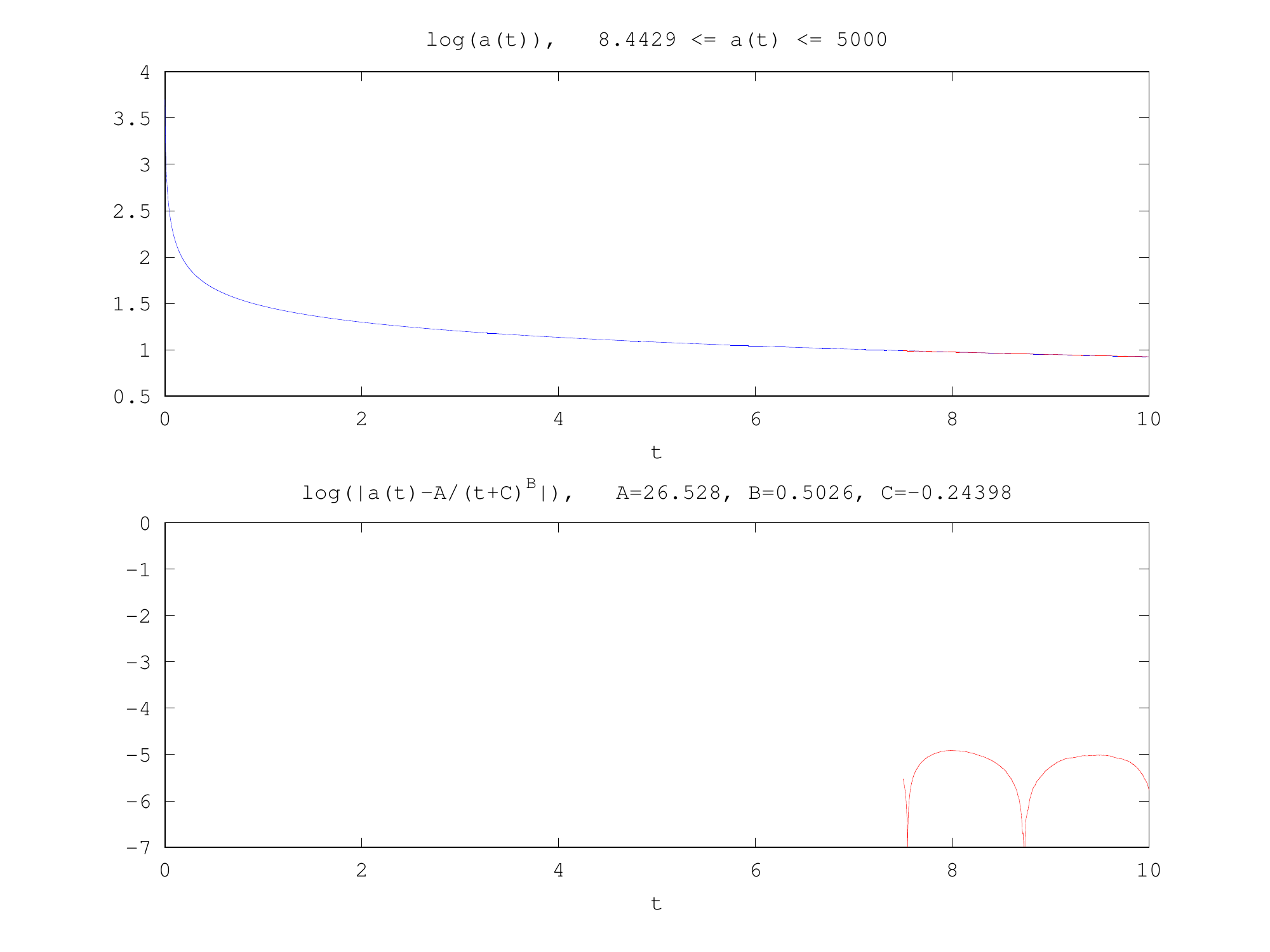}
	\caption{Logarithms of the number density $a(t)$ together with its fit (above), and the fitting error (below)}
	\label{fig1}       
\end{figure}

Figure~\ref{fig1} shows the results on interval [0;10] for $\Delta t=1/3200$. Above the logarithm of $a(t)$ is shown together with the logarithm of its approximation 
\begin{equation}\label{fit}
	a(t) \approx \dfrac{26.528}{(t-0.24389)^{0.5026}}
\end{equation}
on interval [7.5;10]. Below the approximation error less than cca. $10^{-5}$ is shown.

Table~\ref{tab1} compares the results at selected points $t$ for the ``increasing precision'' for the decreasing values $\Delta t$ from 0.01 to 0.0003125. Even for the step $\Delta t=1/3200=0.0003125$ the results may not be considered to have the sufficient precision. If the step number is increasing two times, the calculation time is increasing four times. 

\begin{table}
	\caption{Comparison of the a(t) values at selected points for different values $\Delta t$}
	\label{tab1}
	\centering	
	\begin{tabular}{|l|r|r|r|r|}
		\hline	
		\multicolumn{1}{|c|}{$\Delta t$} & \multicolumn{1}{c|}{$t=0.01$} & \multicolumn{1}{c|}{$t=0.1$} & \multicolumn{1}{c|}{$t=1$} & \multicolumn{1}{c|}{$t=10$} \\ \hline
		0.01      & 2028.8975 & 130.40166 & 41.991715 & 12.607961 \\ \hline
		0.005     & 1338.5228 & 158.18202 & 35.781151 & 10.549856 \\ \hline	
		0.0025    & 1077.5108 & 157.67282 & 32.302138 & \phantom{1}9.379055\\ \hline
		0.00125   & 1062.3410 & 155.70709 & 30.581382 & \phantom{1}8.792275 \\ \hline	 
		0.000625  & 1067.8149 & 154.62549 & 29.841981 & \phantom{1}8.538087 \\ \hline
		0.0003125 & 1068.6433 & 154.18457 & 29.566128 & \phantom{1}8.442899 \\ \hline	
	\end{tabular}
\end{table}

If we suppose, that the difference between the numerical value $a_{\Delta t}(t)$ for the step $\Delta t$ and the ``exact'' value $a(t)$ has the order $p$, e.g., 
\begin{equation}\label{preclow}
	\left|a_{\Delta t}(t) - a(t)\right| \approx C\cdot (\Delta t)^p,
\end{equation}
where both $C$ and $p$ depends on $t$, then using the successive approximations for different steps $\Delta t$, $\Delta t/2$, and $\Delta t/4$ we get
\begin{equation}\label{prec}
	p(t) \approx \log_2\left|\dfrac{a_{\Delta t}(t)-a_{\Delta t/2}(t)}{a_{\Delta t/2}(t)-a_{\Delta t/4}(t)}\right|. 
\end{equation} 

\begin{figure}
	\centering
	\includegraphics[width=0.8\textwidth]{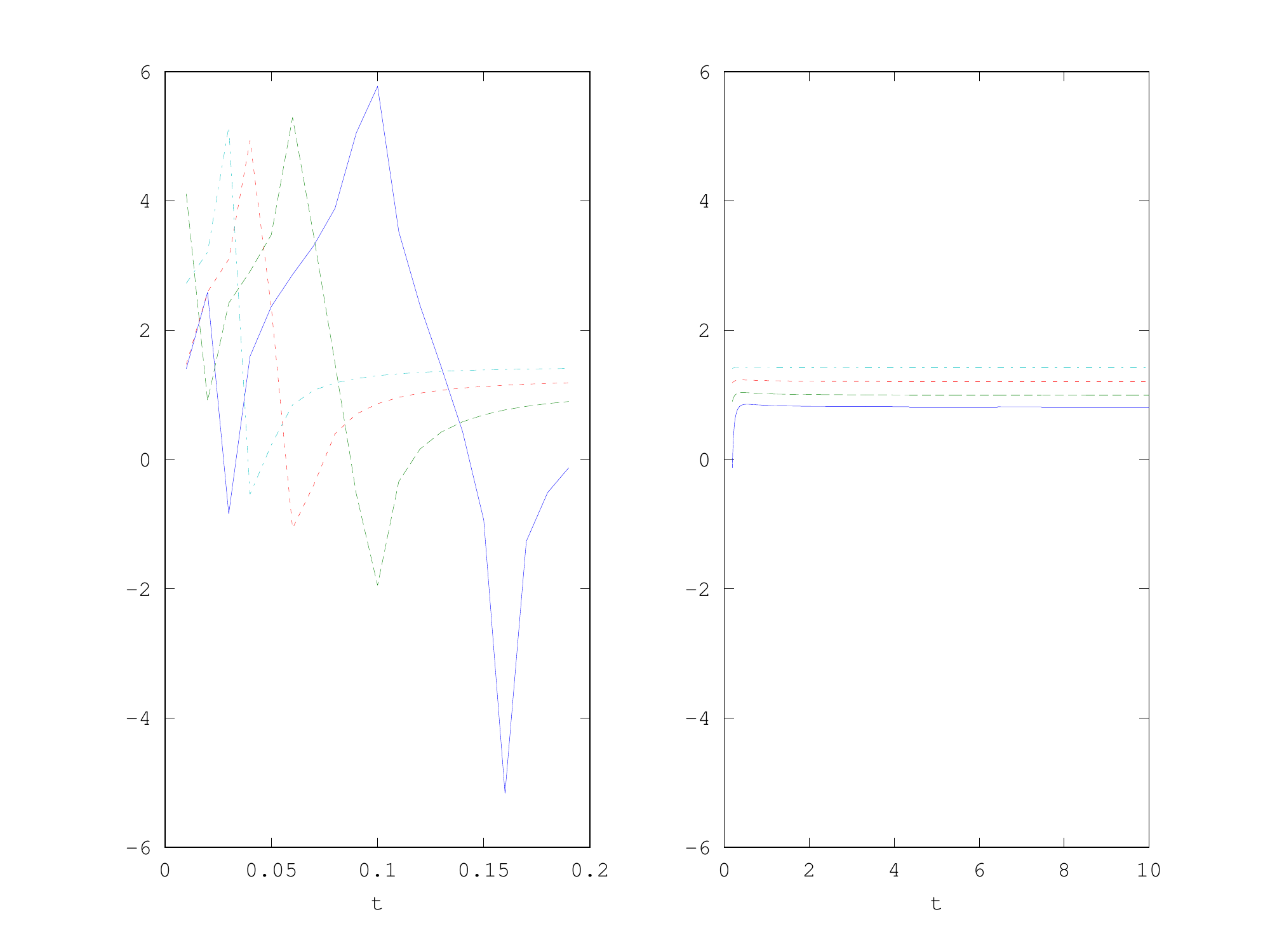}
	\caption{Precision order $p(t)$ for four cases: 0.01--0.005--0.0025, \dots, 0.00125--0.000625--0.0003125}
	\label{fig2}       
\end{figure}

Figure~\ref{fig2} shows the ``precision order'' $p(t)$ calculated for four grid sequences: 0.01--0.005--0.0025, 0.005--0.0025--0.00125, 0.0025--0.00125--0.000625, and 0.00125--0.000625--0.0003125. One can see, that for $t>0.5$ the $p(t)$ behavior is smooth, and the orders are cca. 0.8, 1.0, 1.2, and 1.4, respectively.  
For smaller $t$-values the $p(t)$ behavior is more complicated. 


It seems to be reasonable to study analytically the behavior of the function $a(t)$ for the small values of $t$, and start the numerical computation from some point $t>0$. Also non-uniform grid could be considered.  

\section{Conclusions}
\label{sec-2}

Numerical results presented above could give us some basic imagination about the behavior of the number density function $a(t)$. However, further improvements of the algorithm are necessary.  

It will be also interesting to try to solve the renormalized integro-differential equation 
\[
\frac{\md a(t)}{\md t}=-2\lambda u\nu \mu^{-2\Delta}\, a^2(t)+2\lambda u \nu \mu^{-2\Delta} a^2(t)
\left\{\dfrac{\lambda}{4\pi}\left[\gamma+\ln(2u\nu\mu^2t)\right]\right\}
+
\]
\begin{equation}\label{eqHHL74}
+\dfrac{\lambda^2 u \nu \mu^{-2\Delta}}{2\pi} \int_0^t \frac{\big[a^2(t')-a^2(t)\big]\,\md t'}{t-t'},
\end{equation}
where $\Delta=(d-2)/2$ and $\gamma\doteq 0.57721$ is Euler’s constant, presented in \cite{HHL}, and compare the results.

Another possibility is to study a behavior of the solution of the problem (\ref{eqHHL72}) for $d\to2^-$.

\begin{acknowledgement}
	The work was supported by VEGA Grant 1/0222/13 of the
	Ministry of Education, Science, Research and Sport of the
	Slovak Republic.
\end{acknowledgement}

%
%
%

\end{document}